\documentclass[11pt,onecolumn]{article}
\usepackage{amscd}
\usepackage{times}
\usepackage{amsmath}
\usepackage{amssymb}
\usepackage{xspace}
\usepackage{theorem}
\usepackage{graphicx}
\usepackage{ifpdf}
\usepackage{url,hyperref}
\usepackage{latexsym}
\usepackage{euscript}
\usepackage{xspace}
\usepackage[all]{xy}
\usepackage{color}
\usepackage{makeidx}
\usepackage{tikz}
\long\def\remove#1{}
\newtheorem{theorem}{Theorem}[section] 

\newtheorem{proposition}[theorem]{Proposition}

\newcommand {\mm}[1] {\ifmmode{#1}\else{\mbox{\(#1\)}}\fi}

\newcommand{\img}{\mathrm img}

\newcommand{\cancel}[1]

\begin{document}
 
\title {Comparison between "Parametrized homology via zigzag persistence" 
and "Refined  homology in the presence of a real-valued continuous function"}

\author{
Dan Burghelea  \thanks{
Department of Mathematics,
The Ohio State University, Columbus, OH 43210,USA.
Email: {\tt burghele@math.ohio-state.edu}}
}
\date{}

\date{}
\maketitle
 
\begin{abstract}
\noindent One shows that  the persistence  diagrams  $Dgm^{\cdots} (H\mathbb X)$  and the measures $\mu^{\cdots}_{H\mathbb X}$ defined in \cite {CSKM} can be derived  as restrictions  of the  maps $\delta^f,$ $\gamma^f$ and of the measures $\dim \mathbb F(\cdots),$ $\dim \mathbb T(\cdots)$ considered  in \cite {BH} and \cite{B} to the points and the rectangles above and below diagonal.
\end{abstract}
 
\setcounter{tocdepth}{1}
\tableofcontents

\section {Introduction}

In a recent paper \cite{CSKM} published in this journal  the authors 
have proposed for an $\mathbb R-$ space $\mathbb X$ (= a continuous map  
$f:X\to \mathbb R$), under  reasonable topological hypotheses,  four persistence diagrams 
$$Dgm^{\wedge} (H\mathbb X),  Dgm^{\lor}(H\mathbb X), Dgm^{\backslash \backslash}(H_r \mathbb X), Dgm^{//} (H \mathbb X)$$  regarded as the relevant invariants for parametrized homology of 
$\mathbb X.$  These persistence diagrams  collect the four types of barcodes which can be associated to $f$ via zigzag persistence,  
are maps with discrete support from  $\overline R^2_{+}= \{-\infty \leq  a <b \leq +\infty\}$ to $\mathbb Z_{\geq 0}$  and can be interpreted as  densities $d\mu$ of  four integer valued measures $\mu= \mu^{\wedge}_{H_r(\mathbb X)},$ $ \mu^{\lor}_{H_r(\mathbb X)}, \mu^{\backslash \backslash}_{H_r(\mathbb X)},\mu^{//}_{H_r(\mathbb X)}$ for squares $R= [a,b]\times [c.d],$  $ -\infty\leq a <b <c ,d\leq \infty.$ 

In cite \cite{BH}  sections 6 and 7 (cf. also  \cite {B1}published in this journal) 
and with more details in the book \cite {B} sections 5 and 6, 
under essentially the same hypotheses, two maps with discrete support, $\rho^f_r :\mathbb R^2 \to \mathbb Z_{\geq 0},  \gamma^f_r :\mathbb R^2\setminus \Delta \to \mathbb Z_{\geq 0}$  with 
$\Delta=\{(x,x)\in \mathbb R^2\}$, have been defined and studied.  \footnote { In the above references the discussion refers mostly to  $X$ compact but the conclusions remain the same without the need of any additional hypotheses when $X$ is locally compact and $f: X\to \mathbb R$ is a proper map. In this case the support of the maps $\delta_r^f$ and $\gamma^f_r$ is discrete; this is already used in the case of angle-valued maps}. In the case $X$ is compact the support of the maps $\delta^f_r$ and $\gamma ^f_r$ is finite, hence these maps  are 
configurations of points with multiplicity, cf. \cite{BH}, \cite {B} and \cite{B1}. 

These maps  are derived  from  the vector space-valued maps $\hat \rho^f _r$ and $\hat \gamma^f_r$ with 
$\rho^f= \dim \hat\gamma^f_r$ and $\gamma^f_r= \dim \hat\gamma^f_r$  and  are viewed as relevant topological invariants for a real-valued map $f.$  
The support of $\rho_r$ consists of points in $\mathbb R^2$ which correspond to the closed $r-$bar codes and open $(r-1)-$ barcodes while the support of $\gamma_r$
of points which correspond to the closed-open  and open-closed $r-$barcodes.
In both \cite {BH}  and \cite{B}  these two maps \footnote{They  become four when one treats separately the sigma algebras generated by the above diagonal squares and the below diagonal squares} are defined in two ways, one being  as  "densities" of the measures $\dim \mathbb F_r$ and $\dim \mathbb T_r $ on the sigma algebras  associated to the same type of squares.  
A  measure theoretic  formulation of  both $\delta^f_r, \gamma^f_r,$ and even more general of $\hat\delta^f_r$ and $\hat \gamma^f_r,$  can be explicitly found in \cite{B} subsection 9.2 and is  implicit in \cite{BH} and  \cite {B} subsections 5.1 and 6.1.

The purpose of this note is to show that the persistence diagrams $Dmg^{\cdots}_{\cdots}$ and of the measures $\mu^{\cdots}_{\cdots}$ are the restrictions  to the points and the rectangles above the diagonal resp. below diagonal after composing with the map $T(x,y)= (y,x)$ of  the maps $\rho^f_{\cdots} $ and the measures $\dim\mathbb F_{\cdots},\dim\mathbb T_{\cdots} $ facts  which might pass unnoticed in view of notational differences. 

Precisely, one has: 
\begin{proposition} \label {P1}\ 
\begin{itemize}
\item   (a):  $Dmg^{\wedge}(H_r(\mathbb X))$ equals $\rho^f_r$ restricted to $\mathbb R^2_+, $ 
\item  (b):  $Dmg^{\lor}(H_r(\mathbb X))$ equals  $\rho^f_{r-1}\cdot T$  restricted to $\mathbb R^2_+,$
\item (c):  $Dmg^{\backslash \backslash}(H_r(\mathbb X))$ equals $\rho^f_r$ restricted to $\mathbb R^2_+,$
\item (d): $Dmg^{//}(H_r(\mathbb X))$ equals $\rho^f_{r-1}\cdot T$  restricted to $\mathbb R^2_+.$
\end{itemize}
\end{proposition}
$(R^2_+= \{a,b)\in \mathbb R^2, a<b\}$) and  
\begin{proposition} \label {P2}\ 
\begin{itemize}
\item (A):  $\mu^\wedge_{H_r} ([a,b]\times [c,d))= \dim \mathbb F_r((a,b]\times[c,d)),$ 
\item (B):  $\mu^{\backslash \backslash}_{H_r} ([a,b]\times [c,d])= \dim \mathbb T_r((a,b]\times(c,d]),$ 
\item (C):  $\mu^{//}_{H_{r-1}} ([a,b]\times [c,d])=\dim \mathbb F_r((c,d]\times[a,b)),$ 
\item (D):  $\mu^{\lor}_{H_r} ([a,b]\times [c,d])= \dim \mathbb F_r((c,d]\times[a,b)).$
\end{itemize}
\end{proposition}

Note  that the stability results  as stated  in \cite {CSKM} follows in a straightforward manner  from  the stability results of \cite {BH}  or  \cite {B}
and the Alexander duality in \cite{CSKM} can be derived without effort from the Poincar\'e-duality in \cite {BH} or \cite {B} in  the same way the Alexander-duality 
can be derived from the Poincar\'e duality.  
\vskip .1in 
(NOTE: 
 { After the posting of the first version of this note   it was brought to my attention that the result on Alexander-duality  as well as most of the arguments in \cite{CSKM} were contained in the thesis of one of the author, Sara Kalisnik, cf. http://www.matknjiz.si/doktorati/2013/Kalisnik-14521-4.pdf,  and posted on arXiv cf. Sara Kalisnik , Alexander Duality for Parametrized Homology, arXiv:1303.1591.})

\section {Definitions}
\subsection{CSKM-definitions}

For simplicity in writing one shortens the  notations  in \cite{CSKM} by replacing the notation :  
$\wedge,$ 
$\lor,$ 
$\backslash \backslash,$ 
$//$ 
 by $c,$ $o,$ $co,$ $oc$ abreviatioins of {\it closed, open, closed-open, open-closed} and $\mu^{\cdots} _{H_r\mathbb X}$ by  $\mu^{\cdots}_r.$  

\noindent Consider 

$\mathcal B^c_r=$ the multi-set of closed $r-$bar codes,

$\mathcal B^o_r=$ the multi-set of open $r-$ bar codes,

$\mathcal B^{c,o}_r=$ the multi-set of closed-open $r-$bar codes,

$\mathcal B^{o,c}_r=$ the multi-set of open-closed  $r-$ bar codes.

\noindent The definitions of barcodes  $\mathcal B^{\cdots}_r$ in  \cite{CSKM}  (called in \cite{CSKM} "intervals"  and / or "decorated pairs") are based on the  initial presentation of  {\it zigzag persistence} introduced by Carlsson, de-Silva,  Morozov in 2009 . A reformulation of these  definitions in terms of "death" and of "observability" is provided in \cite{B} subsection  9.1.1.
\vskip .1in

If for a barcode $I$  one denotes by $l(I)$ resp. $r(I)$ the left end resp. the right end, then a careful reading of the definitions in \cite {CSKM}  shows that for a box $R=[a,b]\times [c,d]$ with $-\infty\leq a <b<c <d\leq \infty$ one has:
\begin{equation}
\begin{aligned}
\mu_{H_r}^c(R) =&\sharp \{ I\in \mathcal B_r^c \  \ \mid   a<   l(I) \leq b, \  c\leq r(I) <d\} \\ 
\mu_{H_r}^{c,o}(R) =&\sharp \{ I\in \mathcal B_r^{c,o} \mid   a<   l(I) \leq b,\  c\leq r(I) <d\}\\  
\mu_{H_r}^o(R) =&\sharp \{ I\in \mathcal B_r^o \  \ \mid   a<   l(I) < b,\  c< r(I) <d\}\\  
\mu_{H_r}^{o,c}(R) =&\sharp \{ I\in \mathcal B_r^{o,c} \mid   a<   l(I) \leq b,\  c\leq  r(I) <d\}
\end{aligned}
\end{equation}

and then 
 
\begin{equation}
\begin{aligned} 
Dmg^c_r(a,b)= &\sharp \{ I\in \mathcal B_r^c \mid   a= l(I)\ r(I) =b\} \\  
Dmg^{c,o}_r(a,b)= &\sharp \{ I\in \mathcal B_r^{c,o} \mid   a= l(I)\  r(I) =b\} \\ 
Dmg^o_r(a,b)= &\sharp \{ I\in \mathcal B_r^o \mid   a= l(I)\  r(I) =b\} \\ 
Dmg^{o,c}_r(a,b)=& \sharp \{ I\in \mathcal B_r^{o,c} \mid   a= l(I)\  r(I) =b\} 
\end{aligned}
\end{equation}

\subsection {BH-definitions}

Denote by : 

$\mathbb I_a(r)=\img (H_r(f^{-1}((-\infty, a]))\to H_r(X)),$

 $\mathbb I^a(r)= \img (H_r (f^{-1}([a,\infty)))\to H_r(X)),$  
 
 $\mathbb F_r(a, b)= \mathbb I_a(r)\cap \mathbb I^b(r),$ 
 $F_r(a,b):= \dim \mathbb F_r(a,b),$ 

 $\mathbb T_r(a,b):= \ker (H_r(f^{-1}((-\infty,a]) \to  H_r(f^{-1}((-\infty,b]))$  when   $a<b,$  

$\mathbb T_r(a,b):= \ker (H_r(f^{-1}([a,\infty)) \to  f^{-1}([b, \infty))$ when $a>b,$   

$T_r(a,b)= \dim \mathbb T_r(a,b).$ 
\vskip .1in


\noindent For a box  $B= (a, b]\times [c,d),\ a <b, c<d$ one defines 
$$\mathbb F_r(B):= \mathbb F_r(b,c) /  (\mathbb F_r(a,c)+ \mathbb F_r(b,d));$$
and one observes 
\begin{equation}
\dim \mathbb F_r(B)= F_r(b,c) + F_r(a,d)- F_r(a,c) - F_r(b,d).
\end{equation}
\noindent For a box above diagonal  $B' = (a,b]\times (c,d], a<b \leq c <d$ one defines  
$$\mathbb T_r(B'):= \mathbb T_r({b,d}) / j'  \mathbb T_r{(a,d)}+ \mathbb T_r{(b,c)} ,$$
with $j':  \mathbb T_r{(a,d)}  \to \mathbb T_r{(b,d)}$  
the obviously induced 
linear map, and one observes
\begin{equation}
\dim \mathbb T_r(B')= T_r(b,d) + T_r(a,c)- T_r(a,d) - T_r(b,c).
\end{equation}
\noindent For a box below  diagonal $B''= [c,d)\times [a,b), a<b\leq c <d$  one defines 
 $$\mathbb T_r(B''):= \mathbb T_r(c,a) / j''  \mathbb T_r(d,a)+ \mathbb T_r(c,b)$$
with $j'':  \mathbb T_r({d,a}) \to \mathbb T_r(c,a)$ the obviously induced 
linear map, and one observes 
\begin{equation}
\dim \mathbb T_r(B'')= T_r(c,a) + T_r(d,b)- T_r(c,b) - T_r(d,a).
\end{equation}
\vskip .1in

\noindent Recall from \cite {BD11} Theorem 3.2 and Proposition 5.3 
 or from \cite {BH} Proposition 4.1 \footnote {a real valued map can be regarded as an angle valued map}
 or from \cite {B} Proposition 4.3  the  following equalities:  
\begin{itemize}
\item
 \begin{equation*}
 \dim H_r(f^{-1}(-\infty,a]) = \begin{cases} \sharp \{I\in \mathcal B_r^c\mid l(I) \leq a\} \\
 \sharp \{I\in \mathcal B_{r-1}^o\mid I \subset (-\infty, a)\}\\
 \sharp \{I\in \mathcal B_r^{co}\mid  l(I) \leq a < r(I)\}  
 \end{cases}
 \end{equation*}
 and then 
  \begin{equation*}
 \dim \mathbb I_a(r) = \begin{cases} \sharp \{I\in \mathcal B_r^c\mid l(I) \leq a\} \\
 \sharp \{I\in \mathcal B_{r-1}^o\mid I \subset (-\infty, a) \} 
  \end{cases},
 \end{equation*}
\item 
\begin{equation*}
 \dim H_r(f^{-1}([a,\infty)) = \begin{cases} \sharp \{I\in \mathcal B_r^c\mid r(I) \geq a\} \\
 \sharp \{I\in \mathcal B_{r-1}^o\mid I \subset (a, \infty)\}\\
 \sharp \{I\in \mathcal B_r^{oc}\mid  l(I) < a \leq r(I) \}
\   \end{cases}
 \end{equation*}
 and then 
  \begin{equation*}
 \dim \mathbb I^a(r) = \begin{cases} \sharp \{I\in \mathcal B_r^c\mid r(I) \geq a\} \\
 \sharp \{I\in \mathcal B_{r-1}^o\mid I \subset (a, \infty)\}
 \  \end{cases}.
 \end{equation*}
\end{itemize}
 As a consequence  we have
 \begin{enumerate}
 \item for $a >b$
 \begin{equation}\dim \mathbb F_r(a,b)= 
 \begin{cases} \sharp \{ I\in \mathcal B_r^c\mid l(I)\leq a, r(I)\geq b\} \\
 \sharp \{I\in \mathcal B_{r-1}^o\mid b<l(I) < r(I) <a\}
\   \end{cases}
 \end{equation} 
 \begin{equation}
 \delta^f_r(a,b)=  \sharp \{I=[a,b]\in \mathcal B_{r-1}^o \mid  l(I)= b, r(I)=a\}
\end{equation}

\item for $a\leq b$
\begin{equation}\dim \mathbb F_r(a,b)= 
  \sharp \{I\in \mathcal B_r^c\mid  l(I)\leq a \leq b\leq r(I)\}
  \end{equation}
\begin{equation}
 \delta^f_r(a,b)=  \sharp \{I=[a,b]\in \mathcal B_r^c \mid  l(I)= a, r(I)=b\}
\end{equation}

\item 
for $a<b$
 \begin{equation}
 \dim \mathbb T_{a,b}(r) = 
 \sharp \{I\in \mathcal B_r^{co}\mid l(I) \leq a < r(I)\leq b \}
\end{equation}
\begin{equation}
 \gamma^f_r(a,b)=  \sharp \{I=[a,b)\in \mathcal B_{r-1}^o \mid  l(I)= a, r(I)=b\}
\end{equation}

\item for $a >b$
\begin{equation}
 \dim \mathbb T^{a,b}(r) = 
 \sharp \{I\in \mathcal B_r^{oc}\mid 
 b\leq l(I)< a\leq r(I) 
 \} 
 \end{equation}
\begin{equation}
 \gamma^f_r(a,b)=  \sharp \{I=[a,b]\in \mathcal B_{r-1}^o \mid  l(I)= b, r(I)=a\}
\end{equation}

 \end{enumerate}

\section  {Equalities}

For $R= [a,b]\times [c,d]$ with $a<b< c<d$ 

 (1)+  (8) +(3) imply Proposition\ref{P2} (A)  and  (2)+ (9) imply Proposition\ref{P1}(a),
 
 (1)+  (6) +(3) imply Proposition \ref{P2}(B)  and   (2)+  (7) imply Proposition\ref{P1}(b),
 
 (1)+  (10) +(4) imply Proposition\ref{P2}(C) and (2)+ (11) imply Proposition\ref{P1}(c),
 
 (1)+  (12) +(4) imply Proposition\ref{P2}(D) and (2)+ (13) imply Proposition\ref{P1}(d).

\end{document}